\documentclass[11pt]{article}
\usepackage[pdftex]{hyperref}
\usepackage{amsmath}
\usepackage{mathpazo}
\usepackage{textcomp}
\usepackage{fullpage}
\begin{document}

\begin{center}
\textbf{\large A COMPARATIVE REVIEW OF RECENT RESEARCHES IN GEOMETRY.\footnote[1]{Translated by Dr. M.\ W.\ HASKELL, Assistant Professor of Mathematics in the University of California. Published in Bull. New York Math. Soc. 2, (1892-1893), 215-249. LaTeXed by Nitin C.\ Rughoonauth}}\\
\vspace{0.5cm}
\textnormal{(\it PROGRAMME ON ENTERING THE PHILOSOPHICAL FACULTY AND THE SENATE OF THE UNIVERSITY OF ERLANGEN IN 1872.)}\\
\vspace{0.3in}
BY PROF. FELIX KLEIN.
\end{center}

\indent {\it Prefatory Note by the Author.} - My 1872 Programme, appearing as a separate publication (Erlangen, A.\ Deichert), had but a limited circulation at first. With this I could be satisfied more easily, as the views developed in the Programme could not be expected at first to receive much attention. But now that the general development of mathematics has taken, in the meanwhile, the direction corresponding precisely to these views, and particularly since {\it Lie} has begun the publication in extended form of his {\it Theorie der Transformationsgruppen} (\cite{Li5}), it seems proper to give a wider circulation to the expositions in my Programme. An Italian translation by M. Gina Fano was recently published in the {\it Annali di Matematica}, ser. 2, vol. 17. A kind reception for the English translation, for which I am much indebted to Mr. Haskell, is likewise desired.\\
\indent The translation is an absolutely literal one; in the two or three places where a few words are changed, the new phrases are enclosed in square brackets [ ]. In the same way are indicated a number of additional footnotes which it seemed desirable to append, most of them having already appeared in the Italian translation. - F.\ KLEIN.

MSC-Class: 01A55; 01A75; 53-03; 58-03.

\newpage

Among the advances of the last fifty years in the field of geometry, the development of {\it projective geometry}\footnote[2]{See Note I of the appendix.} occupies the first place. Although it seemed at first as if the so-called metrical relations were not accessible to this treatment, as they do not remain unchanged by projection, we have nevertheless learned recently to regard them also from the projective point of view, so that the projective method now embraces the whole of geometry. But metrical properties are then to be regarded no longer as characteristics of the geometrical figures {\it per se}, but as their relations to a fundamental configuration, the imaginary circle at infinity common to all spheres.

\par When we compare the conception of geometrical figures gradually obtained in this way with the notions of ordinary (elementary) geometry, we are led to look for a general principle in accordance with which the development of both methods has been possible. This question seems the more important as, beside the elementary and the projective geometry, are arrayed a series of other methods, which albeit they are less developed, must be allowed the same right to an individual existence. Such are the geometry of reciprocal radii vectores, the geometry of rational transformations, etc., which will be mentioned and described further on.

\par In undertaking in the following pages to establish such a principle, we shall hardly develop an essentially new idea, but rather formulate clearly what has already been more or less definitely conceived by many others. But it has seemed the more justifiable to publish connective observations of this kind, because geometry, which is after all one in substance, has been only too much broken up in the course of its recent rapid development into a series of almost distinct theories\footnote[3]{See Note II.}, which are advancing in comparative independence of each other. At the same time I was influenced especially by a wish to present certain methods and views that have been developed in recent investigation by {\it Lie} and myself. Our respective investigations, different as has been the nature of the subjects treated, have led to the same generalized conception here presented; so that it has become a sort of necessity to thoroughly discuss this view and on this basis to characterize the contents and general scope of those investigations.

\par Though we have spoken so far only of geometrical investigations, we will include investigations on manifoldnesses of any number of dimensions\footnote[4]{See Note IV.}, which have been developed from geometry by making abstraction from the geometric image, which is not essential for purely mathematical investigations\footnote[5]{See Note III.}. In the investigation of manifoldnesses the same different types occur as in geometry; and, as in geometry, the problem is to bring out what is common and what is distinctive in investigations undertaken independently of each other. Abstractly speaking, it would in what follows be sufficient to speak throughout of manifoldnesses of $n$ dimensions simply; but it will render the exposition simpler and more intelligible to make use of the more familiar space-perceptions. In proceeding from the consideration of geometric objects and developing the general ideas by using these as an example, we follow the path which our science has taken in its development and which it is generally best to pursue in its presentation.

\par A preliminary exposition of the contents of the following pages is here scarcely possible, as it can hardly be presented in a more concise form\footnote[6]{This very conciseness is a defect in the following presentation which I fear will render the understanding of it essentially more difficult. But the difficulty could hardly be removed except by a very much fuller exposition, in which the separate theories, here only touched upon, would have been developed at length.}; the headings of the sections will indicate the general course of thought. 

\par At the end I have added a series of notes, in which I have either developed further single points, wherever the general exposition of the text would seem to demand it, or have tried to define with reference to related points of view the abstract mathematical one predominant in the observations of the text.

\section{\small GROUPS OF SPACE-TRANSFORMATIONS. PRINCIPAL GROUP. FORMULATION OF A GENERAL PROBLEM.}

The most essential idea required in the following discussion is that of a {\it group} of space-trans\-for\-ma\-tions.

\par The combination of any number of tranformations of space\footnote[7]{We always regard the totality of configurations in space as simultaneously affected by the transformations, and speak therefore of {\it transformations of space}. The transformations may introduce other elements in place of points, like dualistic transformations, for instance; there is no distinction in the text in this regard.} is always equivalent to a single transformation. If now a given system of transformations has the property that any transformation obtained by combining any tranformations of the system belongs to that system, it shall be called a {\it group of transformations}\footnote[8]{[This definition is not quite complete, for it has been tacitly assumed that the groups mentioned always include the inverse of every operation they contain; but, when the number of operations is infinite, this is by no means a necessary consequence of the group idea, and this assumption of ours should therefore be explicitly added to the definition of this idea given in the text.] \par The ideas, as well as the notation, are taken from the {\it theory of substitutions}, with the difference merely that there instead of the transformations of a continuous region the permutations of a finite number of discrete quantities are considered.}.

\par An example of a group of transformations is afforded by the totality of motions, every motion being regarded as an operation performed on the whole of space. A group contained in this group is formed, say, by the rotations about one point\footnote[9]{{\it Camille Jordan} has formed all the groups contained in the general group of motions.\cite{Jo}} On the other hand, a group containing the group of motions is presented by the totality of the collineations. But the totality of the dualistic transformations does not form a group; for the combination of two dualistic transformations is equivalent to a collineation. A group is, however, formed by adding the totality of the dualistic to that of the collinear transformations\footnote[10]{It is not at all necessary for the transformations of a group to form a continuous succession, although the groups to be mentioned in the text will indeed always have that property. For example, a group is formed by the finite series of motions which superpose a regular body upon itself, or by the infinite but discrete series which superpose a sine-curve upon itself.}. 

\par Now there are space-transformations by which the geometric properties of configurations in space remain entirely unchanged. For geometric properties are, from their very idea, independent of the position occupied in space by the configuration in question, of its absolute magnitude, and finally of the sense\footnote[11]{By ``sense" is to be understood that peculiarity of the arrangement of the parts of a figure which distinguishes it from the symmetrical figure (the reflected image). Thus, for example, a right-handed and a left-handed helix are of opposite ``sense".} in which its parts are arranged. The properties of a configuration remain therefore unchanged by any motions of space, by transformation into similar configurations, by transformation into symmetrical configurations with regard to a plane (reflection), as well as by any combination of these transformations. The totality of all these transformations we designate as the {\it principal group}\footnote[12]{The fact that these tranformations form a group results from their very idea.} of space-transformations; {\it geometric properties are not changed by the transformations of the principal group}. And, conversely, {\it geometric properties are characterized by their remaining invariant under the transformations of the principal group}. For, if we regard space for the moment as immovable, etc., as a rigid manifoldness, then every figure has an individual character; of all the properties possessed by it as an individual, only the properly geometric ones are preserved in the transformations of the principal group. The idea, here formulated somewhat indefinitely, will be brought out more clearly in the course of the exposition.

\par Let us now dispose with the concrete conception of space, which for the mathematician is not essential, and regard it only as a manifoldness of $n$ dimensions, that is to say, of three dimensions, if we hold to the usual idea of the point as space element. By analogy with the transformations of space we speak of transformations of the manifoldness; they also form groups. But there is no longer, as there is in space, one group distinguished above the rest by its signification; each group is of equal importance with every other. As a generalization of geometry arises then the following comprehensive problem:

\par {\it Given a manifoldness and a group of transformations of the same; to investigate the configurations belonging to the manifoldness with regard to such properties as are not altered by the transformations of the group.}

\par To make use of a modern form of expression, which to be sure is ordinarily used only with reference to a particular group, the group of all the linear transformation, the problem might be stated as follows:

\par {\it Given a manifoldness and a group of transformations of the same; to develop the theory of invariants relating to that group.}

\par This is the general problem, and it comprehends not alone ordinary geometry, but also and in particular the more recent geometrical theories which we propose to discuss, and the different methods of treating manifoldnesses of $n$ dimensions. Particular stress is to laid upon the fact that the choice of the group of transformations to be adjoined is quite arbitrary, and that consequently all the methods of treatment satisfying our general condition are in this sense of equal value.

\section{\small GROUPS OF TRANSFORMATIONS, ONE OF WHICH INCLUDES THE OTHER, ARE SUCCESIVELY ADJOINED. THE DIFFERENT TYPES OF GEOMETRICAL INVESTIGATION AND THEIR RELATION TO EACH OTHER.}

As the geometrical properties of configurations in space remain unaltered under {\it all} the transformations of the principal group, it is by the nature of the question absurd to inquire for such properties as would remain unaltered under only a part of those transformations. This inquiry becomes justified, however, as soon as we investigate the configurations of space in their relation to elements regarded as fixed. Let us, for instance, consider the configurations of space with reference to one particular point, as in spherical trigonometry. The problem then is to develop the properties remaining invariant under the transformations of the principal group, not for the configurations taken independently, but for the system consisting of these configurations together with the given point. But we can state this problem in this other form: to examine configurations in space with regard to such properties as remain unchanged by those transformations of the principal group which can still take place when the point is kept fixed. In other words, it is exactly the same thing whether we investigate the configurations of space taken in connection with the given point from the point of view of the principal group or whether, without any such connection, we replace the principal group by that partial group whose transformations leave the point in question unchanged.

\par This is a principle which we shall frequently apply; we will therefore at once formulate it generally, as follows:

\par Given a manifoldness and a group of transformations applying to it. Let it be proposed to examine the configurations contained in the manifoldness with reference to a given configuration. {\it We may, then, either add the given configuration to the system, and then we have to investigate the properties of the extended system from the point of view of the given group, or we may leave the system unextended, limiting the transformations to be employed to such transformations of the given group as leave the given configuration unchanged. (These transformations necessarily form a group by themselves.)}

\par Let us now consider the converse of the problem proposed at the beginning of this section. This is intelligible from the outset. We inquire what properties of the configurations of space remain unaltered by a group of transformations which contains the principal group as a part of itself. Every property found by an investigation of this kind is a geometric property of the configuration itself; but the converse is not true. In the converse problem we must apply the principle just enunciated, the principal group being now the smaller. We have then:

\par {\it If the principal group be replaced by a more comprehensive group, a part only of the geometric properties remain unchanged. The remainder no longer appear as properties of the configurations of space by themselves, but as properties of the system formed by adding to them some particular configuration.} This latter is defined, in so far as it is a definite\footnote[13]{Such a configuration can be generated, for instance, by applying the transformations of the principal group to any arbitrary element which cannot be converted into itself by any transformation of the given group.} configuration at all, by the following condition: {\it The assumption that it is fixed must restrict us to those transformations of the given group which belong to the principal group.}

\par In this theorem is to be found the peculiarity of the recent geometrical methods to be discussed here, and their relation to the elementary method. What characterizes them is just this, that they base their investigations upon an extended group of space-transformations instead of upon the principal group. Their relation to each other is defined, when one of the groups includes the other, by a corresponding theorem. The same is true of the various methods of treating manifoldnesses of $n$ dimensions which we shall take up. We shall now consider the separate methods from this point of view, and this will afford an opportunity to explain on concrete examples the theorems enunciated in a general form in this and the preceding sections. 

\section{\small PROJECTIVE GEOMETRY.}

Every space-transformation not belonging to the principal group can be used to transfer the properties of known configurations to new ones. Thus we apply the results of plane geometry to the geometry of surfaces that can be represented ({\it abgebildet}) upon a plane;  in this way long before the origin of a true projective geometry the properties of figures derived by projection from a given figure were inferred from those of the given figure. But projective geometry only arose as it became customary to regard the original figure as essentially identical with all those deducible from it by projection, and to enunciate the properties transferred in the process of projection in such a way as to put in evidence their independence of the change due to the projection. By this process {\it the group of all the projective transformations} was made the basis of the theory in the sense of \textsection 1, and that is just what created the antithesis between projective and ordinary geometry. 

\par A course of development similar to the one here described can be regarded as possible in the case of every kind of space-transformation; we shall often refer to it again. It has gone on still further in two directions within the domain of projective geometry itself. On the one hand, the conception was broadened by admitting the {\it dualistic} transformations into the group of the fundamental transformation. From the modern point of view two reciprocal figures are not to be regarded as two distinct figures, but as essentially one and the same. A further advance consisted in extending the fundamental group of collinear and dualistic transformations by the admission in each case of the {\it imaginary} transformations. This step requires that the field of true space-elements has previously been extended so as to include imaginary elements, - just exactly as the admission of dualistic transformations into the fundamental group requires the simultaneous introduction of point and line as space-elements. This is not the place to point out the utility of introducing imaginary elements, by means of which alone we can attain an exact correspondence of the theory of space with the established system of algebraic operations. But, on the other hand, it must be remembered that the reason for introducing the imaginary elements is to be found in the consideration of algebraic operations and not in the group of projective and dualistic transformations. For, just as we can in the latter case limit ourselves to real transformations, since the real collineations and dualistic transformations form a group by themselves, so we can equally well introduce imaginary space-elements even when we are not employing the projective point of view, and indeed must do so in strictly algebraic investigations.

\par How metric properties are to be regarded from the projective point of view is determined by the general theorem of the preceding section. Metrical properties are to be considered as projective relations to a fundamental configuration, the circle at infinity\footnote[14]{This view is to be regarded as one of the most brilliant achievements of [the French school]; for it is precisely what provides a sound foundation for that distinction between properties of position and metrical properties, which furnishes a most desirable starting-point for projective geometry.}, a configuration having the property that it is transformed into itself only by those transformations of the projective group which belong at the same time to the principal group. The proposition thus broadly stated needs a material modification owing to the limitation of the ordinary view taken of geometry as treating only of {\it real} space-elements (and allowing only {\it real} transformations). In order to conform to this point of view, it is necessary expressly to adjoin to the circle at infinity the system of real space-elements (points); properties in the sense of elementary geometry are projectively either properties of the configurations by themselves, or relations to this system of the real elements, or to the circle at infinity, or finally to both.

\par We might here make mention further of the way in which {\it von Staudt}\cite{St1} develops projective geometry, - i.e., that projective geometry which is based on the group containing all the real projective and dualistic transformations\footnote[15]{The extended horizon, which includes {\it imaginary} transformations, was first used by {\it von Staudt} as the basis of his investigation in his later work\cite{St2}.}.

\par We know how, in his system, he selects from the ordinary matter of geometry only such features as are preserved in projective transformations. Should we desire to proceed to the consideration of metrical properties also, what we should have to do would be precisely to introduce these latter as relations to the circle at infinity. The course of thought thus brought to completion is in so far of great importance for the present considerations, as a corresponding development of geometry is possible for every one of the methods we shall take up.

\section{\small TRANSFER OF PROPERTIES BY REPRESENTATIOINS (ABBILDUNG).}

Before going further in the discussion of the geometrical methods which present themselves beside the elementary and the projective geometry, let us develop in a general form certain considerations which will continually recur in the course of the work, and for which a sufficient number of examples are already furnished by the subjects touched upon up to this point. The present section and the following one will be devoted to these discussions.

\par Suppose a manifoldness $A$ has been investigated with reference to a group $B$. If, by any transformation whatever, $A$ be then converted into a second manifoldness $A'$, the group $B$ of transformations, which transformed $A$ into itself, will become a group $B'$, whose transformations are performed upon $A'$. It is then a self-evident principle that {\it the method of treating $A$ with reference to $B$ at once furnishes the method of treating $A'$ with reference to $B'$}, i.e., every property of a configuration contained in $A$ obtained by means of the group $B$ furnishes a property of the corresponding configuration in $A'$ to be obtained by the group $B'$.

\par For example, let $A$ be a straight line and $B$ the $\infty^3$ linear transformations which transform $A$ into itself. The method of treating $A$ is then just what modern algebra designates as the theory of binary forms. Now, we can establish a correspondence between the straight line and a conic section $A'$ in the same plane by projection from a point of the latter. The linear transformations $B$ of the straight line into itself will then become, as can easily be shown, linear transformations $B'$ of the conic into itself, i.e., the changes of the conic resulting from those linear transformations of the plane which transform the conic into itself.

\par Now, by the principle stated in \textsection 2\footnote[16]{The principle might be said to be applied here in a somewhat extended form.}, the study of the geometry of the conic section is the same, whether the conic be regarded as fixed and only those linear transformations of the plane which transform the conic into itself be taken into account, or whether all the linear transformations of the plane be considered and the conic be allowed to vary too. The properties which we recognized in systems of points on the conic are accordingly projective properties in the ordinary sense. Combining this consideration with the result just deduced, we have, then:

\par {\it The theory of binary forms and the projective geometry of systems of points on a conic are one and the same, i.e., to every proposition concerning binary forms corresponds a proposition concerning such systems of points, and vice versa.}\footnote[17]{Instead of the plane conic we may equally well introduce a twisted cubic, or indeed a corresponding configuration in an $n$-dimensional manifoldness.}

\par Another suitable example to illustrate these considerations is the following. If a quadric surface be brought into correspondence with a plane by stereographic projection, the surface will have one fundamental point, - the centre of projection. In the plane there are two, - the projections of the generators passing through the centre of projection. It then follows directly: the linear transformations of the plane which leave the two fundamental points unaltered are converted by the representation ({\it Abbildung}) into linear transformations of the quadric itself, but only into those which leave the centre of projection unaltered. By linear transformations of the surface into itself are here meant the changes undergone by the surface when linear space-transformations are performed which transform the surface into itself. According to this, the projective investigation of a plane with reference to two of its points is identical with the projective investigation of a quadric surface with reference to one of its points. Now, if imaginary elements are also taken into account, the former is nothing else but the investigation of the plane from the point of view of elementary geometry. For the principal group of plane transformations comprises precisely those linear transformations which leave two points (the circular points at infinity) unchanged. We obtain then finally:

\par {\it Elementary plane geometry and the projective investigation of a quadric surface with reference to one of its points are one and the same.}

\par These examples may be multiplied at pleasure\footnote[18]{For other examples, and particularly for the extension to higher dimensions of which those here presented are capable, let me refer to an article of mine\cite{Kl1}, and further to {\it Lie}'s investigations cited later.}; the two here developed were chosen because we shall have occasion to refer to them again.

\section{\small ON THE ARBITRARINESS IN THE CHOICE OF THE SPACE-ELEMENT. HESSE'S PRINCIPLE OF TRANSFERENCE. LINE GEOMETRY.}

As element of the straight line, of the plane, of space, or of any manifoldness to be investigated, we may use instead of the point any configuration contained in the manifoldness, - a group of points, a curve or surface\footnote[19]{See Note III.}, etc. As there is nothing at all determined at the outset about the number of arbitrary parameters upon which these configurations shall depend, the number of dimensioins of our line, plane, space, etc., may be anything we like, according to our choice of the element. {\it But as long as we base our geometrical investigation upon the same group of transformations, the substance of the geometry remains unchanged}. That is to say, every proposition resulting from {\it one} choice of the space-element will be a true proposition under any other assumption; but the arrangement and correlation of the propositions will be changed.

\par The essential thing is, then, the group of transformations; the number of dimensions to be assigned to a manifoldness appears of secondary importance.

\par The combination of this remark with the principle of the last section furnishes many interesting applications, some of which we will now develop, as these examples seem better fitted to explain the meaning of the general theory than any lengthy exposition.

\par Projective geometry on the straight line (the theory of binary forms) is, by the last section, equivalent to projective geometry on the conic. Let us now regard as element on the conic the point-pair instead of the point. Now, the totality of the point-pairs of the conic may be brought into correspondence with the totality of the straight lines in the plane, by letting every line correspond to that point-pair in which it intersects the conic. By this representation ({\it Abbildung}) the linear transformations of the conic into itself are converted into those linear transformations of the plane (regarded as made up of straight lines) which leave the conic unaltered. But whether we consider the group of the latter, or whether we base our investigation on the totality of the linear transformations of the plane, always adjoining the conic to the plane configurations under investigation, is by \textsection 2 one and the same thing. Uniting all these considerations, we have:

\par {\it The theory of binary forms and projective geometry of the plane with reference to a conic are identical}.

\par Finally, as projective geometry of the plane with reference to a conic, by reason of the equality of its group, coincides with that projective metrical geometry which in the plane can be based upon a conic\footnote[20]{See Note V.}, we can also say:

\par {\it The theory of binary forms and general projective metrical geometry in the plane are one and the same}.

\par In the preceding consideration the conic in the plane might be replaced by the twisted cubic, etc., but we will not carry this out further. The correlation here explained between the geometry of the plane, of space, or of a manifoldness of any number of dimensions is essentially identical with the principle of transference proposed by {\it Hesse} (Borchardt's Journal, vol. 66).

\par An example of much the same kind is furnished by the projective geometry of space; or, in other words, the theory of quaternary forms. If the straight line be taken as space-element and be determined, as in line geometry, by six homogeneous co-ordinates connected by a quadratic equation of condition, the linear and dualistic transformations of space are seen to be those linear transformations of the six variables (regarded as independent) which transform the equation of condition into itself. By a combination of considerations similar to those just developed, we obtain the following theorem:

\par {\it The theory of quaternary forms is equivalent to projective measurement in a manifoldness generated by six homogeneous variables}. 

\par For a detailed exposition of this view I will refer to \cite{Kl2}, and to a note at the close of this paper\footnote[21]{See Note VI.}.

\par To the foregoing expositions I will append two remarks, the first of which is to be sure implicitly contained in what has already been said, but needs to be brought out at length, because the subject to which it applies is only too likely to be misunderstood.

\par Through the introduction of arbitrary configurations as space-elements, space becomes of any number of dimensions we like. But if we then keep to the (elementary or projective) space-perception with which we are familiar, the fundamental group for the manifoldness of $n$ dimensions is given at the outset; in the one case it is the principal group, in the other the group of projective transformations. If we wished to take a different group as a basis, we should have to depart from the ordinary (or from the projective) space-perception. Thus, while it is correct to say that, with a proper choice of space-elements, space represents manifoldnesses of any number of dimensions, it is equally important to add that {\it in this representation either a definite group must form the basis of the investigation of the manifoldness, or else, if we wish to choose the group, we must broaden our geometrical perception accordingly}. If this were overlooked, an interpretation of line geometry, for instance, might be sought in the following way. In line geometry the straight line has six co-ordinates: the conic in the plane has the same number of coefficients. The interpretation of line geometry would then be the geometry in a system of conics separated from the aggregation of all conics by a quadratic equation between the coefficients. This is correct, provided we take as fundamental group for the plane geometry the totality of the transformations represented by the linear transformations of the coefficients of the conic which transform the quadratic equation into itself. But if we retain the elementary or the projective view of plane geometry, we have no interpretation at all.

\par The second remark has reference to the following line of reasoning: Suppose in space some group or other, the principal group for instance, be given. Let us then select a single configuration, say a point, or a straight line, or even an ellipsoid, etc., and apply to it all the transformations of the principal group. We thus obtain an infinite manifoldness with a number of dimensions in general equal to the number of arbitrary parameters contained in the group, but reducing in special cases, namely, when the configuration originally selected has the property of being transformed into itself by an infinite number of the transformations of the group. Every manifoldness generated in this way may be called, with reference to the generating group, a {\it body}\footnote[22]{In choosing this name I follow the precedent established by {\it Dedekind} in the theory of numbers, where he applies the name {\it body} to a system of numbers formed from given elements by given operations (\cite{Di})}.

\par If now we desire to base our investigations upon the group, selecting at the same time certain definite configurations as space-elements, and if we wish to represent uniformly things which are of like characteristics, {\it we must evidently choose our space-elements in such a way that their manifoldness either is itself a body or can be decomposed into bodies}. This remark, whose correctness is evident, will find application later (\textsection 9). This idea of a body will come under discussion once more in the closing section, in connection with certain related ideas\footnote[23]{[In the text sufficient attention is not paid to the fact that the proposed group may contain so-called self-conjugate subgroups. If a geometrical configuration remain unchanged by the operations of a self-conjugate subgroup, the same is true for all configurations into which it is transformed by the operations of the whole group; i.e., for all configurations of the body arising from it. But a body so formed would be absolutely unsuited to represent the operations of the group. In the text, therefore, are to be admitted only bodies formed of space-elements which remain unchanged by no self-conjugate subgroup of the given group whatever.]}.

\section{\small THE GEOMETRY OF RECIPROCAL RADII. INTERPRETATION OF $x+iy$.}

With this section we return to the discussion of the various lines of geometric research, which was begun in \textsection\textsection 2 and 3.

\par As a parallel in many respects to the processes of projective geometry, we may consider a class of geometric investigations in which the transformation by reciprocal radii vectores (geometric inversion) is continually employed. To these belong investigations on the so-called eyelides and other anallagmatic surfaces, on the general theory of orthogonal systems, likewise on potential, etc. It is true that the processes here involved have not yet, like projective geometry, been united into a special geometry, {\it whose fundamental group would be the totality of the transformations resulting from a combination of the principal group with geometric inversion}; but this may be ascribed to the fact that the theories named have never happened to receive a connected treatment. To the individual investigators in this line of work some such systematic conception can hardly have been foreign.

\par The parallel between this geometry of reciprocal radii and projective geometry is apparent as soon as the question is raised; it will therefore be sufficient to call attention in a general way to the following points:

\par In projective geometry the elementary ideas are the point, line, and plane. The circle and the sphere are but special cases of the conic section and the quadric surface. The region at infinity of elementary geometry appears as a plane; the fundamental configuration to which elementary geometry is referred is an imaginary conic at infinity.

\par In the geometry of reciprocal radii the elementary ideas are the point, circle, and sphere. The line and the plane are special cases of the latter, characterized by the property that they contain a point which, however, has no further special significance in the theory, namely, the point at infinity. If we regard this point as fixed, elementary geometry is the result.

\par The geometry of reciprocal radii admits of being stated in a form which places it alongside of the theory of binary forms and of line geometry, provided the latter be treated in the way indicated in the last section. To this end we will for the present restrict our observations to plane geometry and therefore to the geometry of reciprocal radii in the plane\footnote[24]{The geometry of reciprocal radii on the straight line is equivalent to the projective investigation of the line, as the transformations in question are the same. Thus in the geometry of reciprocal radii, also, we can speak of the anharmonic ratio of four points on a line and of four points on a circle.}.

\par We have already referred to the connection between elementary plane geometry and the projective geometry of the quadric surface with one distinctive point (\textsection 4). If we disregard the distinctive point, that is to say, if we consider the projective geometry on the surface by itself, we have a representation of the geometry of reciprocal radii in the plane. For it is easy to see\footnote[25]{See \cite{Kl1}.} that to the group of geometric inversion in the plane corresponds by virtue of the representation ({\it Abbildung}) of the quadric surface the totality of the linear transformations of the latter into itself. We have, therefore,

\par {\it The geometry of reciprocal radii in the place and the projective geometry on a quadric surface are one and the same}; and, similarly:

\par {\it The geometry of reciprocal radii in space is equivalent to the projective treatment of a manifoldness represented by a quadratic equation between five homogeneous variables}.

\par By means of the geometry of reciprocal radii space geometry is thus brought into exactly the same connection with a manifoldness of four dimensions as by means of [projective] geometry with a manifoldness of five dimensions.

\par The geometry of reciprocal radii in the plane, if we limit ourselves to {\it real} transformations, admits of an interesting interpretation, or application, in still another direction. For, representing the complex variable $x+iy$ in the plane in the usual way, to its linear transformations corresponds the group of geometric inversion, with the above-mentioned restriction to real operations\footnote[26]{[The language of the text is inexact. To the linear transformations $z'=\frac{\alpha z+\beta}{\gamma z+\delta}$ (where $z'=x'+iy'$, $z=x+iy$) correspond only those operations of the group of geometric inversion by which no reversion of the angles takes place (in which the two circular points of the plane are not interchanged). If we wish to include the whole group of geometric inversion, we must, in addition to the transformations mentioned, take account of the other (not less important) ones given by the formula $z'=\frac{\alpha \bar{z}+\beta}{\gamma \bar{z}+\delta}$ (where again $z'=x'+iy'$, but $\bar{z}=x-iy$).]}. But the investigation of functions of a complex variable, regarded as subject to any linear transformations whatever, is merely what, under a somewhat different mode of representation, is called the theory of binary forms. In other words:

\par {\it The theory of binary forms finds interpretation in the geometry of reciprocal radii in the real plane, and precisely in the way in which complex values of the variables are represented}.

\par From the plane we will ascend to the quadric surface, to return to the more familiar circle of ideas of the projective transformations. As we have taken into consideration only real elements of the plane, it is not a matter of indifference how the surface is chosen; it can evidently not be a ruled surface. In particular, we may regard it as a spherical surface, - as is customary for the interpretation of a complex variable, - and obtain in this way the theorem:

\par {\it The theory of the binary forms of a complex variable finds representation in the projective geometry of the real spherical surface}.

\par I could not refrain from setting forth in a note\footnote[27]{See Note VII.} how admirably this interpretation illustrates the theory of binary cubics and quartics.

\section{\small EXTENSION OF THE PRECEDING CONSIDERATIONS. LIE'S SPHERE GEOMETRY.} 

With the theory of binary forms, the geometry of reciprocal radii, and line geometry, which in the foregoing pages appear co-ordinated and only distinguished by the number of variables, may be connected certain further developments, which shall now be explained. In the first place, these developments are intended to illustrate with new exmaples the idea that the group determining the treatment of given subjects can be extended indefinitely; but, in the second place, the intention was particularly to explain the relation to the views here set forth of certain considerations presented by {\it Lie} in a recent article\cite{Li2}. The way by which we here arrive at {\it Lie}'s sphere geometry differs in this respect from the one pursued by {\it Lie}, that he proceeds from the conceptions of line geometry, while we assume a smaller number of variables in our exposition. This will enable us to be in agreement with the usual geometric perception and to preserve the connection with what precedes. The investigation is independent of the number of variables, as {\it Lie} himself has already pointed out (\cite{Li1a}\cite{Li1b}). It belongs to that great class of investigations concerned with the projective discussion of quadratic equations between any number of variables,  - investigations upon which we have already touched several times, and which will repeatedly meet us again (see \textsection 10, for instance).

\par I proceed from the connection established between the real plane and the sphere by stereographic projection. In \textsection 5 we connected plane geometry with the geometry on a conic section by making the straight line in the plane correspond to the point-pair in which it meets the conic. Similarly we can establish a connection between space geometry and the geometry on the sphere, by letting every plane of space correspond to the circle in which it cuts the sphere. If then by stereographic projection we transfer the geometry on the sphere from the latter to the plane (every circle being thereby transformed into a circle), we have the following correspondence:

\par the space geometry whose element is the plane and whose group is formed of the linear transformations converting a sphere into itself, and

\par the plane geometry whose element is the circle and whose group is the group of geometric inversion.

\par The former geometry we will now generalize in two directions by substituting for its group a more comprehensive group. The resulting extension may then be immediately transferred to plane geometry by representation ({\it Abbildung}).

\par Instead of those linear transformations of space (regarded as made up of planes) which convert the sphere into itself, it readily suggests itself to select either the totality of the {\it linear} transformations of space, or the totality of those plane-transformations which leave the sphere unchanged [in a sense yet to be examined]; in the former case we dispense with the sphere, in the latter with the linear character of the transformations. The former generalization is intelligible without further explanation; we will therefore consider it first and follow out its importance for plane geometry. To the second case we shall return later, and shall then in the first place have to determine the most general transformation of that kind.

\par Linear space-transformations have the common property of converting pencils and sheafs of planes into like pencils and sheafs. Now, transferred to the sphere, the pencil of planes gives a pencil of circles, i.e., a system of $\infty^1$ circles with common intersections; the sheaf of planes gives a sheaf of circles, i.e., a system of $\infty^2$ circles perpendicular to a fixed circle (the circle whose plane is the polar plane of the point common to the planes of the given sheaf). Hence to linear space-transformations there correspond on the sphere, and furthermore in the plane, circle-transformations characterized by the property that they convert pencils and sheafs of circles into the same\footnote[29]{Such transformations are considered in Grassmann\cite{Gr1b}.} {\it The plane geometry which employs the group of transformations thus obtained is the representation of ordinary projective space geometry}. In this geometry the point cannot be used as element of the plane, for the points do not form a {\it body} (\textsection 5) for the chosen group of transformations; but circles shall be chosen as elements.

\par In the case of the second extension named, the first question to be settled is with regard to the nature of the group of transformations in question. The problem is, to find plane-transformations converting every [pencil] of planes whose [axis touches] the sphere into a like [pencil]. For brevity of expression, we will first consider the reciprocal problem and, moreover, go down a step in the number of dimensions; we will therefore look for point-transformations of the plane which convert every tangent to a given conic into a like tangent. To this end we regard the plane with its conic as the representation of a quadric surface projected on the plane from a point of space not in the surface in such a way that the conic in question represents the boundary curve. To the tangents to the conic correspond the generators of the surface, and the problem is reduced to that of finding the totality of the point-transformations of the surface into itself by which generators remain generators.

\par Now, the number of these transformations is, to be sure, $\infty^n$, where $n$ may have any value. For we only need to regard the point on the surface as intersection of the generators of the two systems, and to transform each system of lines into itself in any way whatever. But among these are in particular the linear transformations, and to these alone will we attend. For, if we had to do, not with a surface, but with an $n$-dimensional manifoldness represented by a quadratic equation, the linear transformations alone would remain, the rest would disappear\footnote[30]{If the manifoldness be stereographically projected, we obtain the well-known theorem: in regions of $n$ dimensions (even in space) there are no isogonal point-transformations except the transformations of the group of geometric inversion. In the plane, on the other hand, there are any number besides. See the articles by {\it Lie} already cited.}.

\par These linear transformations of the surface into itself, transferred to the plane by projection (other than stereographic), give two-valued point-transformations, by which from every tangent to the boundary conic is produced, it is true, a tangent, but from every other straight line in general a conic having double contact with the boundary curve. This group of transformations will be conveniently characterized by basing a projective measurement upon the boundary conic. The transformations will then have the property of converting points whose distance apart is zero by this measurement, and also points whose distance from a given point is constant, into points having the same properties. 

\par All these considerations may be extended to any number of variables, and can in particular be applied to the original inquiry, which had reference to the sphere and plane as elements. We can then give the result an especially perspicuous form, because the angle formed by two planes according to the projective measurement referred to a sphere is equal to the angle in the ordinary sense formed by the circles in which they intersect the sphere.

\par We thus obtain upon the sphere, and furthermore in the plane, a group of circle-transformations having the property that {\it they convert circles which are tangent to each other (include a zero angle), and also circles making equal angles with another circle, into like circles}. The group of these transformations contains on the sphere the linear transformations, in the plane the transformations of the group of geometric inversion\footnote[31]{[Perhaps the addition of some few analytic formulae will materially help to explain the remarks in the text. Let the equation of the sphere, which we project stereographically on the plane, be in ordinary tetrahedral co-ordinates:
\begin{equation}
x_1^2+x_2^2+x_3^2+x_4^2=0.\nonumber
\end{equation}
The $x$'s satisfying this equation of condition we then interpret as tetracyclic co-ordinates in the plane.
\begin{equation}
u_1 x_1+u_2 x_2+u_3 x_3+u_4 x_4=0\nonumber
\end{equation}
will be the general circular equation of the plane. If we compute the radius of the circle represented in this way, we come upon the square root $\sqrt{u_1^2+u_2^2+u_3^2+u_4^2}$, which we may denote by $i u_3$. We can now regard the circles as elements of the plane. The group of geometric inversion is then represented by the totality of those homogeneous linear transformations of $u_1, u_2, u_3, u_4$, by which $u_1^2+u_2^2+u_3^2+u_4^2$ is converted into a multiple of itself. But the extended group which corresponds to {\it Lie}'s sphere geometry consists of those homogeneous linear transformations of the five variables $u_1, u_2, u_3, u_4$, which convert $u_1^2+u_2^2+u_3^2+u_4^2+u_5^2$ into a multiple of itself.]}.

\par The circle geometry based on this group is analogous to the sphere geometry which {\it Lie} has devised for space and which appears of particular importance for investigations on the curvature of surfaces. It includes the geometry of reciprocal radii in the same sense as the latter includes elementary geometry.

\par The circle- (sphere-) transformations thus obtained have, in particular, the property of converting circles (spheres) which touch each other into circles (spheres) having the same property. If we regard all curves (surfaces) as envelopes of circles (spheres), then it results from this fact that curves (surfaces) which touch each other will always be transformed into curves (surfaces) having the same property. The transformations in question belong, therefore, to the class of {\it contact-transformations} to be considered from a general standpoint further on, i.e., transformations under which the contact of point-configurations is an invariant relation. The first circle-transformations mentioned in the present section, which find their parallel in corresponding sphere-transformations, are not contact transformations.

\par While these two kinds of generalization have here been applied only to the geometry of reciprocal redii, they nevertheless hold in a similar way for line geometry and in general for the projective investigation of a manifoldness defined by a quadratic equation, as we have already indicated, but shall not develop further in this connection.

\section{\small ENUMERATION OF OTHER METHODS BASED ON A GROUP OF\\ POINT-TRANSFORMATIONS.}

Elementary geometry, the geometry of reciprocal radii, and likewise projective geometry, if we disregard the dualistic transformations connected with the interchange of the space-element, are included as special cases among the large number of conceivable methods based on groups of point-transformations. We will here mention especially only the three following methods, which agree in this respect with those named. Though these methods are far from having been developed into independent theories in the same degree as projective geometry, yet they can clearly be traced in the more recent investigations\footnote[32]{[Groups with a finite number of parameters having been treated in the examples hitherto taken up, the so-called infinite groups will now be the subject of consideration in the text.]}.

\subsection{The Group of Rational Transformations.} 

\par In the case of rational transformations we must carefully distinguish whether they are rational for {\it all} points of the region under consideration, viz., of space, or of the plane, etc., or only for the points of a manifoldness contained in the region, viz., a surface or curve. The former alone are to be employed when the problem is to develop a geometry of space or of the plane in the meaning hitherto understood; the latter obtain a meaning, from our point of view, only when we wish to study the geometry on a given surface or curve. The same distinction is to be drawn in the case of the {\it analysis situs} to be discussed presently.

\par The investigations in both subjects up to this time have been occupied mainly with transformations of the second kind. Since in these investigations the question has not been with regard to the geometry on the surface or curve, but rather to find the criteria for the transformability of two surfaces or curves into each other, they are to be excluded from the sphere of the investigations here to be considered\footnote[33]{[From another point of view they are brought back again, which I did not yet know in 1872, very nicely into connection with the considerations in the text. Given any algebraic configuration (curve, or surface, etc.), let it be transferred into a higher space by introducing the ratios
\begin{equation}
\phi_1 : \phi_2 : \ldots : \phi_P \nonumber
\end{equation}
of the intergrands of the first species belonging to it as homogeneous co-ordinates. In this space we have then simply to take the group of homogeneous linear transformations as a basis for our further considerations. See various articles by {\it Brill}, {\it N\"{o}ther}, and {\it Weber}, and (to mention a single recent article) my own paper: \cite{Kl3}.]}. For the general synopsis here outlined does not embrace the entire field of mathematical research, but only brings certain lines of thought under a common point of view.

\par Of such a geometry of rational transformations as must result on the basis of the transformations of the first kind, only a beginning has so far been made. In the region of the first grade, viz., on the straight line, the rational transformations are identical with the linear transformations and therefore furnish nothing new. In the plane we know the totality of rational transformations (the Cremona transformations); we know that they can be produced by a combination of quadratic transformations. We know further certain invariant properties of plane curves [with reference to the totality of rational transformations], viz., their deficiency, the existence of moduli; but these considerations have not yet been developed into a geometry of the plane, properly speaking, in the meaning here intended. In space the whole theory is still in its infancy. We know at present but few of the rational transformations, and use them to establish correspondences between known and unknown surfaces.

\subsection{Analysis situs.}

\par In the so-called analysis situs we try to find what remains unchanged under transformations resulting from a combination of infinitesimal distortions. Here, again, we must distinguish whether the whole region, all space, for instance, is to be subjected to the transformations, or only a manifoldness contained in the same, a surface. It is the transformations of the first kind on which we could found a space geometry. Their group would be entirely different in constitution from the groups heretofore considered. Embracing as it does all transformations compounded from (real) infinitesimal point-transformations, it necessarily involves the limitation to real space-elements, and belongs to the domain of arbitrary functions. This group of transformations can be extended to advantage by combining it with those real collineations which at the same time affect the region at infinity.

\subsection{The Group of all Point-transformations.}

\par While with reference to this group no surface possesses any individual characteristics, as any surface can be converted into any other by transformations of the group, the group can be employed to advantage in the investigation of higher configurations. Under the view of geometry upon which we have taken our stand, it is a matter of no importance that these configurations have hitherto not been regarded as geometric, but only as analytic, configurations, admitting occasionally of geometric application, and, furthermore, that in their investigation methods have been employed (these very point-transformations, for instance) which we have only recently begun to consciously regard as geometric transformations. To these analytic configurations belong, above all, homogeneous differential expressions, and also partial differential equations. For the general discussion of the latter, however, as will be explained in detail in the next section, the more comprehensive group of all contact-transformations seems to be more advantageous.

\par The principal theorem in force in the geometry founded on the group of all point-trans\-for\-ma\-tions is this: {\it that for an infinitesimal portion of space a point-trans\-for\-ma\-tion always has the value of a linear transformation}. Thus the developments of projective geometry will have their meaning for infinitesimals; and, whatever be the choice of the group for the treatment of the manifoldness, {\it in this fact lies a distinguishing characteristic of the projective view}.

\par Not having spoken for some time of the relation of methods of treatment founded on groups, one of which includes the other, let us now give one more example of the general theory of \textsection 2. We will consider the question how projective properties are to be understood from the point of view of ``all point-transformations," disregarding here the dualistic transformations which, properly speaking, form part of the group of projective geometry. This question is identical with the other question, What condition differentiates the group of linear point-transformations from the totality of point-transformations? What characterizes the linear group is this, that to every plane it makes correspond a plane; it contains those transformations under which the manifoldness of planes (or, what amounts to the same thing, of straight lines) remains unchanged. {\it Projective geometry is to be obtained from the geometry of all point-transformations by adjoining the manifoldness of planes, just as elementary is obtained from projective geometry by adjoining the imaginary circle at infinity}. Thus, for instance, from the point of view of all point-transformations the designation of a surface as an algebraic surface of a certain order must be regarded as an invariant relation to the manifoldness of planes. This becomes very clear if we connect, as {\it Grassmann} does\cite{Gr2}, the generation of algebraic configurations with their construction by lines.

\section{\small ON THE GROUP OF ALL CONTACT-TRANSFORMATIONS.}

Particular cases of contact-transformations have been long known; {\it Jacobi} has even made use of the most general contact-transformations in analytical investigations, but an effective geometrical interpretation has only been given them by recent researches of {\it Lie}'s\footnote[34]{See, in particular, the article already cited\cite{Li2}. For the details given in the text in regard to partial differential equations I am indebted mainly to oral communications of {\it Lie}'s; see his note, {\it Zur Theorie partieller Differentialgleichungen}\cite{Li3}.}. It will therefore not be superfluous to explain here in detail what a contact-transformation is. In this we restrict ourselves, as hitherto, to point-space with its three dimensions. 

\par By a contact-transformation is to be understood, analytically speaking, any substitution which expresses the values of the variables $x, y, z$ and their partial derivatives $\frac{dz}{dx}=p$, $\frac{dz}{dy}=q$ in terms of new variables $x', y', z', p', q'$. It is evident that such substitutions, in general, convert surfaces that are in contact into surfaces in contact, and this accounts for the name. Contact-transformations are divided into three classes (the point being taken as space-element), viz., those in which {\it points} correspond to the $\infty^3$ points (the point-transformations just considered); those converting the points into curves; lastly, those converting them into surfaces. This classification is not to be regarded as essential, inasmuch as for other $\infty^3$ space-elements, say for planes, while a division into three classes again occurs, it does not coincide with the division occurring under the assumption of points as elements.

\par If a point be subjected to all contact-transformations it is converted into the totality of points, curves, and surfaces. Only in their entirety, then, do points, curves, and surfaces form a {\it body} of our group. From this may be deduced the general rule that the formal treatment of a problem from the point of view of all contact-transformations (e.g., the theory of partial differential equations considered below) must be incomplete if we operate only with point- (or plane-) co-ordinates, for the very reason that the chosen space-elements do not form a body.

\par If, however, we wish to preserve the connection with the ordinary methods, it will not do to introduce as space-elements all the individual configurations contained in the body, as their number is $\infty^\infty$. This makes it necessary to introduce in these considerations as space-element not the point, curve, or surface, but the ``surface-element," i.e., the system of values $x, y, z, p, q$. Each contact-transformation converts every surface-element into another; the $\infty^5$ surface-elements accordingly form a body.

\par From this point of view, point, curve, and surface must be uniformly regarded as aggregates of surface-elements, and indeed of $\infty^2$ elements. For the surface is covered by $\infty^2$ elements, the curve is tangent to the same number, through the point pass the same number. But these aggregates of $\infty^2$ elements have another characteristic property in common. Let us designate as the {\it united position} of two consecutive surface-elements $x, y, z, p, q$ and $x+dx, y+dy, z+dz, p+dp, q+dq$ the relation defined by the equation
\begin{equation}
dz-pdx-qdy=0. \nonumber
\end{equation}
Thus point, curve, and surface agree in being {\it manifoldnesses of $\infty^2$ elements, each of which is united in position with the $\infty^1$ adjoining elements}. This is the common characteristic of point, curve, and surface; and this must serve as the basis of the analytical investigation, if the group of contact-transformations is to be used.

\par The united position of consecutive elements is an invariant relation under any contact-trans\-for\-ma\-tion whatever. And, conversely, contact-transformations may be defined as {\it those substitutions of the five variables $x, y, z, p, q$, by which the relation
\begin{equation}
dz-pdx-qdy=0 \nonumber
\end{equation}
is converted into itself}. In these investigations space is therefore to be regarded as a manifoldness of five dimensions; and this manifoldness is to be treated by taking as fundamental group the totality of the transformations of the variables which leave a certain relation between the differentials unaltered.

\par First of all present themselves as subjects of investigation the manifoldnesses defined by one or more equations between the variables, i.e., {\it by partial differential equations of the first order, and systems of such equations}. It will be one of the principal problems to select out of the manifoldnesses of elements satisfying given equation systems of $\infty^1$, or of $\infty^2$, elements which are all united in position with a neighboring element. A question of this kind forms the sum and substance of the problem of the solution of a partial differential equation of the first order. It can be formulated in the following way: to select from among the $\infty^4$ elements satisfying the equation all the twofold manifoldnesses of the given kind. The problem of the complete solution thus assumes the definite form: to classify in some way the $\infty^4$ elements satisfying the equation into $\infty^2$ manifoldnesses of the given kind.

\par It cannot be my intention to pursue this consideration of partial differential equations further; on this point I refer to {\it Lie}'s articles already cited. I will only point out one thing further, that from the point of view of the contact-transformations a partial differential equation of the first order has no invariant, that every such equation can be converted into any other, and that therefore linear equations in particular have no distinctive properties. Distinctions appear only when we return to the point of view of the point-transformations.

\par The groups of contact-transformations, of point-transformations, finally of projective transformations, may be defined in a uniform manner which should here not be passed over\footnote[35]{I am indebted to a remark of {\it Lie}'s for these definitions.}. Contact-transformations have already been defined as those transformations under which the united position of consecutive surface-elements is preserved. But, on the other hand, point-transformations have the characteristic property of converting consecutive line-elements which are united in position into line-elements similarly situated; and, finally, linear and dualistic transformations maintain the united position of consecutive connex-elements. By a connex-element is meant the combination of a surface-element with a line-element contained in it; consecutive connex-elements are said to be united in position when not only the point but also the line-element of one is contained in the surface-element of the other. The term connex-element (though only preliminary) has reference to the configurations recently introduced into geometry by {\it Clebsch}\footnote[36]{\cite{Cl1} and especially \cite{Cl2}.} and represented by an equation containing simultaneously a series of point-coordinates as well as a series of plane- and a series of line-coordinates whose analogues in the plane {\it Clebsch} denotes as connexes.

\section{\small ON MANIFOLDNESSES OF ANY NUMBER OF DIMENSIONS.}

We have already repeatedly laid stress on the fact that in connecting the expositions thus far with space-perception we have only been influenced by the desire to be able to develop abstract ideas more easily through dependence on graphic examples. But the considerations are in their nature independent of the concrete image, and belong to that general field of mathematical research which is designated as the theory of manifoldnesses of any dimensions, - called by {\it Grassmann} briefly ``theory of extension" (Ausdehnungslehre). How the transference of the preceding development from space to the simple idea of a manifoldness is to be accomplished is obvious. It may be mentioned once more in this connection that in the abstract investigation we have the advantage over geometry of being able to choose arbitrarily the fundamental group of transformations, while in geometry a minimum group - the principal group - was given at the outset.

\par We will here touch, and that very briefly, only on the following three methods: 

\subsection{The Projective Method or Modern Algebra (Theory of Invariants).}

Its group consists of the totality of linear and dualistic transformations of the variables employed to represent individual configurations in the manifoldness; it is the generalization of projective geometry. We have already noticed the application of this method in the discussion of infinitesimals in a manifoldness of one more dimension. It includes the two other methods to be mentioned, in so far as its group includes the groups upon which those methods are based.

\subsection{The Manifoldness of Constant Curvature.}

The notion of such a manifoldness arose in {\it Riemann}'s theory from the more general idea of a manifoldness in which a differential expression in the variables is given. In his theory the group consists of the totality of those transformations of the variables which leave the given expression unchanged. On the other hand, the idea of a manifoldness of constant curvature presents itself when a projective measurement is based upon a given quadratic equation between the variables. From this point of view as compared with {\it Riemann}'s the extension arises that the variables are regarded as complex; the variability can be limited to the real domain afterwards. Under this head belong the long series of investigations touched on in \textsection\textsection 5, 6, 7.

\subsection{The Plane Manifoldness.}

{\it Riemann} designates as a plane manifoldness one of constant zero curvature. Its theory is the immediate generalization of elementary geometry. Its group can, like the principal group of geometry, be separated from out the group of the projective method by supposing a configuration to remain fixed which is defined by two equations, a linear and a quadratic equation. We have then to distinguish between real and imaginary if we wish to adhere to the form in which the theory is usually presented. Under this head are to be counted, in the first place, elementary geometry itself, then for instance the recent generalizations of the ordinary theory of curvature, etc.

\begin{center}
\section*{\small CONCLUDING REMARKS.}
\end{center}

In conclusion we will introduce two further remarks closely related to what has thus far been presented, -  one with reference to the analytic form in which the ideas developed in the preceding pages are to be represented, the other marking certain problems whose investigation would appear important and fruitful in the light of the expositions here given.

\par Analytic geometry has often been reproached with giving preference to arbitrary elements by the introduction of the system of co-ordinates, and this objection applied equally well to every method of treating manifoldnesses in which individual configurations are characterized by the values of variables. But while this objection has been too often justified owing to the defective way in which, particularly in former times, the method of co-ordinates was manipulated, yet it disappears when the method is rationally treated. The analytical expressions arising in the investigation of a manifoldness with reference to its group must, from their meaning, be independent of the choice of the co-ordinate system; and the problem is then to clearly set forth this independence analytically. That this can be done, and how it is to be done, is shown by modern algebra, in which the abstract idea of an invariant that we have here in view has reached its clearest expression. It possesses a general and exhaustive law for constructing invariant expressions, and operates only with such expressions. This object should be kept in view in any formal (analytical) treatment, even when other groups than the projective group form the basis of the treatment\footnote[37]{[For instance, in the case of the groups of rotations of three-dimensional space about a fixed point, such a formalism is furnished by quaternions.]}. For the analytical formulation should, after all, be congruent with the conceptions whether it be our purpose to use it only as a precise and perpicuous expression of the conceptions, or to penetrate by its aid into still unexplored regions.

\par The further problems which we wished to mention arise on comparing the views here set forth with the so-called {\it Galois} theory of equations.

\par In the {\it Galois} theory, as in ours, the interest centres on groups of transformations. The objects to which the transformations are applied are indeed different; there we have to do with a finite number of discrete elements, here with the infinite number of elements in a continuous manifoldness. But still the comparison may be pursued further owing to the identity of the group-idea\footnote[38]{I should like here to call to mind {\it Grassmann}'s comparison of combinatory analysis and extensive algebra in the introduction to \cite{Gr1a}.}, an I am the more ready to point it out here, as it will enable us to characterize the position to be awarded to certain investigations begun by {\it Lie} and myself\footnote[39]{See \cite{KL}} in accordance with the views here developed.

\par In the Galois theory, as it is presented for instance in {\it Serret}'s \cite{Se} or in {\it C. Jordan}'s \cite{Jo1}, the real subject of investigation is the group theory of substitution theory itself, from which the theory of equations results as an application. Similarly we require a {\it theory of transformations}, a theory of the groups producible by transformations of any given characteristics. The ideas of commutativity, of similarity, etc., will find application just as in the theory of substitutions. As an application of the theory of transformations appears that treatment of a manifoldness which results from taking as a basis the groups of transformations.

\par In the theory of equations the first subjects to engage the attention are the symmetric functions of the coefficients, and in the next place those expressions which remain unaltered, if not under all, yet under a considerable number of permutations of the roots. In treating a manifoldness on the basis of a group our first inquiry is similarly with regard to the bodies (\textsection 5), viz., the configurations which remain unaltered under all the transformations of the group. But there are configurations admitting not all but some of the transformations of the group, and they are next of particular interest from the point of view of the treatment based on the group; they have distinctive characteristics. It amounts, then, to distinguishing in the sense of ordinary geometry symmetric and regular bodies, surfaces of revolution and helicoidal surfaces. If the subject be regarded from the point of view of projective geometry, and if it be further required that the transformations converting the configurations into themselves shall be commutative, we arrive at the configurations considered by {\it Lie} and myself in the article cited, and the general problem proposed in \textsection 6 of that article. The determination (given in \textsection\textsection 1, 3 of that article) of all groups of an infinite number of commutative linear transformations in the plane forms a part of the general theory of transformations named above\footnote[40]{I must refrain from referring in the text to the fruitfulness of the consideration of infinitesimal transformations in the theory of differential equations. In \textsection 7 of the article cited, {\it Lie} and I have shown that ordinary differential equations which admit the same infinitesimal transformations present like difficulties of integration. How the considerations are to be employed for partial differential equations, {\it Lie} has illustrated by various examples in several places; for instance, in the article named above (\cite{Li2}). See in particular \cite{Li4}.
\par [At this time I may be allowed to refer to the fact that it is exactly the two problems mentioned in the text which have influenced a large part of the further investigations of {\it Lie} and myself. I have already called attention to the appearance of the two first volumes of {\it Lie}'s ``Theorie der Transformationsgruppen."\cite{Li5} Of my own work might be mentioned the later researches on regular bodies, on elliptic modular functions, and on single-valued functions with linear transformations into themselves, in general. An account of the first of these was given in a special work\cite{Kl4}; an exposition of the theory of the elliptic modular functions, elaborated by {\it Dr. Fricke} is in course of publication.]}.\\

\section*{NOTES.}

\subsection*{I. On the Antithesis between the Synthetic and the Analytic Method in Modern Geometry.}

The distinction between modern synthesis and modern analytic geometry must no longer be regarded as essential, inasmuch as both subject-matter and methods of reasoning have gradually taken a similar form in both. We choose therefore in the text as common designation of them both the term {\it projective geometry}. Although the synthetic method has more to do with space-perception and thereby imparts a rare charm to its first simple developments, the realm of space-perception is nevertheless not closed to the analytic method, and the formulae of analytic geometry can be looked upon as a precise and perspicuous statement of geometrical relations. On the other hand, the advantage to original research of a well formulated analysis should not be underestimated, - an advantage due to its moving, so to speak, in advance of the thought. But it should always be insisted that a mathematical subject is not to be considered exhausted until it has become intuitively evident, and the progress made by the aid of analysis is only a first, though a very important, step.

\subsection*{II. Division of Modern Geometry into Theories.}

When we consider, for instance, how persistently the mathematical physicist disregards the advantages afforded him in many cases by only a moderate cultivation of the projective view, and how, on the other hand, the student of projective geometry leaves untouched the rich mine of mathematical truths brought to light by the theory of the curvature of surfaces, we must regard the present state of mathematical knowledge as exceedingly incomplete and, it is to be hoped, as transitory.

\subsection*{III. On the Value of Space-perception.}

When in the text we designated space-perception as something incidental, we meant this with regard to the purely mathematical contents of the ideas to be formulated. Space-perception has then only the value of illustration, which is to be estimated very highly from the pedagogical stand-point, it is true. A geometric model, for instance, is from this point of view very instructive and interesting.

\par But the question of the value of space-perception in itself is quite another matter. I regard it as an independent question. There is a true geometry which is not, like the investigations discussed in the text, intended to be merely an illustrative form of more abstract investigations. Its problem is to grasp the full reality of the figures of space, and to interpret - and this is the mathematical side of the question - the relations holding for them as evident results of the axioms of space-perception. A model, whether constructed and observed or only vividly imagined, is for this geometry not a means to an end, but the subject itself.

\par This presentation of geometry as an independent subject, apart from and independent of pure mathematics, is nothing new, of course. But it is desirable to lay stress explicitly upon this point of view once more, as modern research passes it over almost entirely. This is connected with the fact that, {\it vice versa}, modern research has seldom been employed in investigations on the form-relations of space-configurations, while it appears well adapted to this purpose.

\subsection*{IV. On Manifoldnesses of any Number of Dimensions.}

That space, regarded as the locus of points, has only three dimensions, does not need to be discussed from the mathematical point of view; but just as little can anybody be prevented from that point of view from claiming that space really has four, or any unlimited number of dimensions, and that we are only able to perceive three. The theory of manifoldnesses, advancing as it does with the course of time more and more into the foreground of modern mathematical research, is by its nature fully independent of any suh claim. But a nomenclature has become established in this theory which has indeed been derived from this idea. Instead of the elements of a manifoldness we speak of the points of a higher space, etc. The nomenclature itself has certain advantages, in that it facilitates the interpretation by calling to mind the perceptions of geometry. but it has had the unfortunate result of causing the whide-spread opinion that investigations on manifoldnesses of any number of dimensions are inseparably connected with the above-mentioned idea of the nature of space. Nothing is more unsound than this opinion. The mathematical investigations in question would, it is true, find an immediate application to geometry, if the idea were correct; but their value and purport is absolutely independent of this idea, and depends only on their own mathematical contents.

\par It is quite another matter when {\it Pl\"{u}cker} shows how to regard actual space as a manifoldness of any number of dimensions by introducing as space-element a configuration depending on any number of parameters, a curve, surface, etc. (see \textsection 5 of the text).

\par The conception in which the element of a manifoldness (of any number of dimensions) is regarded as analogous to the point in space was first developed, I suppose, by {\it Grassmann} in his ``Ausdehnungslehre" (\cite{Gr1a}). With him the thought is absolutely free of the above-mentioned idea of the nature of space; this idea goes back to occasional remarks by {\it Gauss}, and became more widely known through {\it Riemann}'s investigations on manifoldnesses, with which it was interwoven.

\par Both conceptions - {\it Grassmann}'s as well as {\it Pl\"{u}cker}'s - have their own peculiar advantages; they can be alternately employed with good results.

\subsection*{V. On the So-called Non-Euclidean Geometry.}

The projective metrical geometry alluded to in the text is essentially coincident, as recent investigations have shown, with the metrical geometry which can be developed under non-acceptance of the axiom of parallels, and is to-day under the name of non-Euclidean geometry widely treated and discussed. The reason why this name has not been mentioned at all in the text, is closely related to the expositions given in the preceding note. With the name non-Euclidean geometry have been associated a multitude of non-mathematical ideas, which have been as zealously cherished by some as resolutely rejected by others, but with which our purely mathematical considerations have nothing to do whatever. A wish to contribute towards clearer ideas in this matter has occasioned the following explanations.

\par The investigations referred to on the theory of parallels, with the results growing out of them, have a definite value for mathematics from two points of views.

\par They show, in the first place, - and this function of theirs may be regarded as concluded once for all, - that the axiom of parallels is not a mathematical consequence of the other axioms usually assumed, but the expression of an essentially new principle of space-perception, which has not been touched upon in the foregoing investigations. Similar investigations could and should be performed with regard to every axiom (and not alone in geometry); an insight would thus be obtained into the mutual relation of the axioms.

\par But, in the second place, these investigations have given us an important mathematical idea, - the idea of a manifoldness of constant curvature. This idea is very intimately connected, as has already been remarked and in \textsection 10 of the text discussed more in detail, with the projective measurement which has arisen independently of any theory of parallels. Not only is the study of this measurement in itself of great mathematical interest, admitting of numerous applications, but it has the additional feature of including the measurement given in geometry as a special (limiting) case and of teaching us how to regard the latter from a broader point of view.

\par Quite independent of the views set forth is the question, what reasons support the axiom on parallels, i.e., whether we should regard it as absolutely given, as some claim, or only as approximately proved by experience, as others say. Should there be reasons for assuming the latter position, the mathematical investigations referred to afford us then immediately the means for constructing a more exact geometry. But the inquiry is evidently a philosophical one and concerns the most general foundations of our understanding. The  mathematician as such is not concerned with this inquiry, and does not wish his investigations to be regarded as dependent on the answer given to the question from the one of the other point of view\footnote[41]{[To the explanations in the text I should like to add here two supplementary remarks.

\par In the first place, when I say that the mathematician as such has no stand to take place on the philosophical question, I do not mean to say that the philosopher can dispense with the mathematical developments in treating the aspect of the question which interests him; on the contrary, it is my decided conviction that a study of these developments is the indispensable prerequisite to every philosophical discussion of the subject.

\par Secondly, I have not meant to say that my {\it personal} interest is exhausted by the mathematical aspect of the question. For my conception of the subject, in general, let me refer to a recent paper\cite{Kl5}]}.

\subsection*{VI. Line Geometry as the Investigation of a Manifoldness of Constant Curvature.}

In combining line geometry with the projective measurement in a manifoldness of five dimensions, we must remember that the straight lines represent elements of the manifoldness which, metrically speaking, are at infinity. It then becomes necessary to consider what the value of a system of a projective measurement is for the elements at infinity; and this may here be set forth somewhat at length, in order to remove any difficulties which might else seem to stand in the way of conceiving of line geometry as a metrical geometry. We shall illustrate these expositions by the graphic example of the projective measurement based on a quadric surface.

\par Any two points in space have with respect to the surface an absolute invariant, - the anharmonic ratio of the two points together with the two points of intersection of the line joining them with the surface. But when the two points move up to the surface, this anharmonic ratio becomes zero independently of the position of the points, except in the case where the two points fall upon a generator, when it becomes indeterminate. This is the only special case which can occur in their relative position unless they coincide, and we have therefore the theorem:

\par {\it The projective measurement in space based upon a quadric surface does not yet furnish a measurement for the geometry on the surface}.

\par This is connected with the fact that by linear transformations of the surface into itself any three points of the surface can be brought into coincidence with three others\footnote[42]{These relations are different in ordinary metrical geometry; for there it is true that two points at infinity have an absolute invariant. The contradiction which might thus be found in the enumeration of the linear transformations of the surface at infinity into itself is removed by the fact that the translations and transformations of similarity contained in this group do not alter the region at infinity at all.}.

\par If a measurement on the surface itself be desired, we must limit the group of transformations, and this result is obtained by supposing any arbitrary point of space (or its polar plane) to be fixed. Let us first take a point not on the surface. We can then project the surface from the point upon a plane, when a conic will appear as the boundary curve. Upon this conic we can base a projective measurement in the plane, which must then be transferred back to the surface\footnote[43]{See \textsection 7 of the text.}. This is a measurement with constant curvature in the true sense, and we have then the theorem: 

\par {\it Such a measurement on the surface is obtained by keeping fixed a point not on the surface}.

\par Correspondingly, we find\footnote[44]{See \textsection 4 of the text}:
\par {\it A measurement with zero curvature on the surface is obtained by choosing as the fixed point a point of the surface itself}.

\par In all these measurements on the surface the generators of the surface are lines of zero length. The expression for the element of arc on the surface differs therefore only by a factor in the different cases. There is no absolute element of arc upon the surface; but we can of course speak of the angle formed by two directions on the surface.

\par All these theorems and considerations can now be applied immediately to line geometry. Line-space itself admits at the outset no measurement, properly speaking. A measurement is only obtained by regarding a linear complex as fixed; and the measurement is of constant or zero curvature, according as the complex is a general or a special one (a line). The selection of a particular complex carries with it further the acceptation of an absolute element of arc. Independently of this, the directions to adjoining lines cutting the given line arc of zero length, and we can besides speak of the angle between any two directions\cite{Kl1}.

\subsection*{VII. On the Interpretation of Binary Forms.}

We shall now consider the graphic illustration which can be given to the theory of invariants of binary cubics and biquadratics by taking advantage of the representation of $x+iy$ on the sphere.

\par A binary cubic $f$ has a cubic covariant $Q$, a quadratic covariant $\Delta$, and an invariant $R$\footnote[46]{See in this connection the corresponding sections of Clebsch\cite{Cl3}}. From $f$ and $Q$ a whole system of covariant sextics $Q^2+\lambda R f^2$ may be compounded, among them being $\Delta^3$. It can be shown\footnote[47]{By considering the linear transformations of $f$ into itself. See \cite{Kl1a}.} that every covariant of the cubic must resolve itself into such groups of six points. Inasmuch as $\lambda$ can assume complex values, the number of these covariants is $\infty^2$\cite{B}.

\par The whole system of forms thus defined can now be represented upon the sphere as follows. By a suitable linear transformation of the sphere into itself let the three points representing $f$ be converted into three equidistant points of a great circle. let this great circle be denoted as the equator, and let the three points $f$ have the longitudes $0^\circ$, $120^\circ$, $240^\circ$. Then $Q$ will be  represented by the points of the equator whose longitudes are $60^\circ$, $180^\circ$, $300^\circ$; $\Delta$ by the two poles. Every form $Q^2+\lambda R f^2$ is represented by six points, whose latitude and longitude are given in the following table, where $\alpha$ and $\beta$ are arbitrary numbers:

\begin{center}
\begin{tabular*}{6in}{|p{1in}|p{1in}|p{1in}|p{1in}|p{1in}|p{1in}|}
\hline
$\alpha$&$\alpha$&$\alpha$&$-\alpha$&$-\alpha$&$-\alpha$\\
$\beta$&$120^\circ+\beta$&$240^\circ+\beta$&$-\beta$&$120^\circ-\beta$&$240^\circ-\beta$\\
\hline
\end{tabular*}
\end{center} 

\par In studying the variation of these systems of points on the sphere, it is interesting to see how they give rise to $f$ and $Q$ (each reckoned twice) and $\Delta$ (reckoned three times).

\par A biquadratic $f$ has a biquadratic covariant $H$, a sextic covariant $T$, and two invariants $i$ and $j$. Particularly noteworthy is the pencil of biquadratic forms $iH+\lambda j f$, all belonging to the same $T$, among them being the three quadratic factors into which $T$ can be resolved, each reckond twice.

\par Let the centre of the sphere now be taken as the origin of a set of rectangular axes $OX$, $OY$, $OZ$. Their six points of intersection with the sphere make up the form $T$. The four points of a set $iH+\lambda j f$ are given by the following table, $x, y, z$ being the co-ordinates of any point of the sphere:

\begin{center}
\begin{tabular}{ c c c c c c c c }
$$ & $x,$ & $\hspace{0.3in}$ & $$ & $y,$ & $\hspace{0.3in}$ & $$ & $z,$\\
$$ & $x,$ & $\hspace{0.3in} $ & $-$ & $y,$ & $\hspace{0.3in}$ & $-$ & $z,$\\
$-$ & $x,$ & $\hspace{0.3in}$ & $$ & $y,$ & $\hspace{0.3in}$ & $-$ & $z,$\\
$-$ & $x,$ & $\hspace{0.3in}$ & $-$ & $y,$ & $\hspace{0.3in}$ & $$ & $z.$\\
\end{tabular}
\end{center}
 
\par The four points are in each case the vertices of a symmetrical tetrahedron, whose opposite edges are bisected by the co-ordinate axes; and this indicates the r$\hat{o}$le played by $T$ in the theory of biquadratic equations as the resolvent of $iH+\lambda j f$.

\vspace{0.5in}

ERLANGEN, {\it October}, 1872.

\end{document}